# A Truncated Manuscript

*Pierre Schapira*

*Récoltes et Semailles I, II. Réflexions et témoignage sur un passé de mathématicien*
by Alexander Grothendieck
Editions Gallimard, 29,50 €.

STRICTLY SPEAKING, this essay is not solely a review of Alexander Grothendieck's *Récoltes et Semailles* (Reaping and Sowing).[1] Although the book as a whole, as well as Grothendieck's work and life, will be discussed here, a good part of the essay is devoted to refuting a thesis Grothendieck developed throughout his original text. We have an important ally in this respect: the author himself. In a series of important additions that have not been incorporated in this new version, Grothendieck goes back entirely on some of his assertions.

THE FIGURE of Grothendieck dominated much of mathematics during the second half of the twentieth century. If his work is essentially concerned with algebraic geometry, his vision and methods have spread far beyond—to algebraic topology, representation theory, complex geometry, symplectic geometry, algebraic analysis and even, more recently, computational geometry. In short, all linear mathematics, as opposed to dynamical systems, probabilities, or purely geometric geometry, such as Riemannian or Hamiltonian geometry. It was under his influence that the language of derived categories and sheaf theory were established in these fields. It was also Grothendieck who had the intuition and then formulated the main lines of the theory of derived categories. He left it to his student Jean-Louis Verdier to write out all the details for his thesis, in which Verdier clarified the key notion of a triangulated category.[2] But it was Grothendieck who situated sheaf theory and the six operations in the theory of derived categories, to which we will return later.

Between 1950 and 1970, functions—possibly generalized—were studied on real or complex manifolds, and especially on Euclidean spaces, using the Fourier transform. But on a complex manifold, when we refer to functions we really mean holomorphic functions, and these present a serious difficulty. That is, they do not exist, at least not globally on a compact manifold such as the projective line—apart, of course, from the constants. Global knowledge, therefore, does not provide any information, unlike the real differentiable case, and one must work locally. There is a rather extraordinary tool for this purpose: sheaf theory, which was invented by Jean Leray while he was a prisoner of war in Germany between 1941 and 1945.[3] Leray's text was somewhat incomprehensible but was later clarified by Henri Cartan and Jean-Pierre Serre, resulting in Roger Godement's famous book, *Théorie des Faisceaux* (Theory of Sheaves).[4] Cartan and Serre used this tool in their respective studies of holomorphic functions in dimension ≥ 1, after the seminal work of Kiyoshi Oka, giving rise to Cartan's Theorems A and B and Serre's influential paper "*Faisceaux algébriques cohérents*" (Coherent algebraic sheaves).[5]

It is in this context that Grothendieck approached algebraic geometry around 1955, providing a solid basis for sheaf cohomology in a foundational paper published by the *Tōhoku Mathematical Journal*.[6] In this paper, one already encounters—implicitly—the main difficulty of category theory, namely the problem of universes, a problem later solved by Grothendieck in *SGA 4*. This was done in the manner of another Alexander cutting the Gordian knot:[7] Grothendieck poses the axiom that every set belongs to a universe. This problem of universes, also known as the inaccessible cardinals, outside of which category theory cannot develop, is probably the reason why Bourbaki gave up on categories and why Grothendieck then left the group. Ralf Krömer has written an excellent article on this issue.[8]

During the 1940s and 1950s, there were two conceptual revolutions whose importance was not immediately understood: sheaf theory, as mentioned above, and category theory, the latter due to Samuel Eilenberg and Saunders Mac Lane.[9] Moreover, the categorical point of view is part of a vast movement of ideas that embraced the structuralist approach of Claude Lévi-Strauss and the linguistics of Noam Chomsky. Instead of considering sets endowed with certain structures, category theory focuses on the relations that can exist between objects. A category $C$ is thus a family of objects—as a set is a family of elements—but given two objects $X$ and $Y$, there exists a priori a set called $\text{Hom}_C(X, Y)$ representing the mor-phisms from *X* to *Y*, these data being, of course, subject to a certain number of natural axioms—composition of morphisms, identity morphisms, etc. A new step then consists in looking at morphisms between categories, which are called functors. Certain key notions then emerge—such as those of adjoint functors, final or initial objects, limits and colimits—giving a precise and unifying meaning to many ideas that run through mathematics.

There is a family of categories that plays a central role: these are the additive categories and, among them, the abelian categories, which are modeled on the category of modules over a ring. But if the vector spaces over a field are replaced by the modules over a ring, the classical tensor product and internal Hom functors are no longer exact, i.e., they do not transform exact sequences into exact sequences—a subspace does not always admit a supplementary. It is thus necessary to consider the derived functors. We then enter the domain of homological algebra, a natural generalization of linear algebra. Here, the reference book was initially that of Cartan and Eilenberg,[10] before it was dethroned by Grothendieck's *Tōhoku* paper. But computing the derived functor of the composition of two functors requires the use of Leray's spectral sequences, giving rise to often inextricable calculations. This is where derived categories show their power: in this language, everything is remarkably simple.

What are the six operations? With ordinary functions, we have three natural operations, besides addition: the product and, associated with an application $f : X \to Y$ between real manifolds, the integration which sends—modulo some technical details—functions on *X* to functions on *Y* and the composition by f which sends functions on *Y* to functions on *X*. In sheaf theory, the tensor product $\overset{L}{\otimes}$ is the analogue of the product, the proper direct image $Rf_!$ is the analogue of the integration, and the inverse image $f^{-1}$ is the analogue of the composition by *f*. But the tensor product has a right adjoint, Rhom, the functor $f^{-1}$ has a right adjoint, the direct image $Rf_*$, and the functor $Rf_!$ has a right adjoint, $f^!$.

The functor $f^!$, which exists only in the derived framework, unlike the other five, was discovered by Grothendieck in the context of étale cohomology and was subsequently constructed by Verdier for locally compact spaces. As Grothendieck had seen, $f^!$ provides a broad generalization of Poincaré duality and this functor now plays a crucial role. But since locally compact spaces appear more often than the étale topology, it is the name Poincaré–Verdier, if not just Verdier alone, that remains associated with duality. This attribution is largely unfair and left Grothendieck feeling somewhat bitter—rightly so.[11]

One might think that such an abstract framework dispenses with explicit computations, but this is a misconception: put simply, the computations are no longer the same. If the direct image functor does not allow for integrals to be computed explicitly, the formalism of the six operations nevertheless gives rise to sophisticated numerical results, such as the computation of dimensions of cohomology spaces. The Riemann–Roch–Hirzebruch–Grothendieck theorem is a beautiful illustration.

In the same vein, one of Grothendieck's fundamental discoveries was to develop sheaf theory on categories and thus, in particular, on spaces that no longer have any points. What do the sheaves require to exist? Namely, the data of open sets and their inclusions, and the notion of a covering. There is nothing to prevent the objects of a category from playing the role of the open sets, the category is then called a pre-site, and it remains to define axiomatically what the coverings are in order to obtain a site—i.e., a category with a Grothendieck topology. This natural generalization of the usual topological spaces proves to be extremely fruitful and analysts would, in fact, do well to draw inspiration from it. On a real manifold, there are far too many pathological open sets and too many coverings if one is interested in what happens at the edge of an open set.

One then arrives at topos theory—topoi for the scholars. The underlying idea, which in a particular situation goes back to Israel Gelfand, is that a space—in this case, a site—can be reconstructed from the category of sheaves on that site. A topos is then a category equivalent to a category of sheaves. The category of sets, for example, is nothing other than the topos associated with a point. But even if topos theory has been used in a new proof of Paul Cohen's result on the independence of the continuum hypothesis, its applications in mathematics are still uncertain.

This presentation of a selection of Grothendieck's fundamental ideas is far from complete and only reflects the particular interests of the reviewer. In *R&S*, Grothendieck lists what he considers the 12 key ideas of his work and also provides a list of his students.[12]

One should, of course, also mention his first works in functional analysis, dating from around 1955, and scheme theory, which revolutionized algebraic geometry, as well as the intuition of motives, a partially conjectural theory that was later developed by Pierre Deligne, Vladimir Voevodsky, Joseph Ayoub, and many others. One should also mention the fundamental text "*A la poursuite des champs*" in which Grothendieck lays the foundations for ∞-categories and homotopical algebra.[13] Indeed, if triangulated categories are an incredibly simple and efficient tool, they have a defect which seriously limits their use, for example, in gluing problems. This defect is linked to the fact that a certain morphism is unique up to isomorphism, but this isomorphism is not unique! The new theory of ∞-categories, to which the names Jacob Lurie, Graeme Segal, Bertrand Töen, and a few others must be associated, is in the process of completely supplanting the classical theory of derived categories, although it is, for the time being, not easily accessible, to say the least.

PIERRE CARTIER has written a remarkable article on Grothendieck's life and it is pointless to paraphrase it here.[14] There are also excellent articles by Allyn

Jackson and Winfried Scharlau, as well as all the links on the Grothendieck Circle website managed by Leila Schneps.[15]

Nonetheless, a few words on this subject are helpful. Grothendieck's father, Sascha Schapiro, was a Russian anarchist who took part in the aborted revolution of 1905. He then served ten years in the prisons of Czar Nicholas II before being released following the revolution of 1917. Despite initially being feted as a hero, Schapiro was soon declared an enemy of the people. He later fought alongside the Republicans during the Spanish Civil war before becoming a traveling photographer in France. In 1939, Schapiro was interned at the Camp Vernet in the French Pyrenees. He was handed over to the Nazis by the Vichy police in 1942 and disappeared into Auschwitz.

Grothendieck's mother, Hanka, was an extreme-left militant in Germany during the 1920s who emigrated to France when Adolf Hitler came to power.[16] Her son did not join her until 1938, at the age of 10, after having lived in hiding on a farm in Germany. Grothendieck spent part of the war in Le Chambon-sur-Lignon at the famous Collège Cévenol that saved so many Jewish children.

Grothendieck's mathematical life began in Nancy during the 1950s, where Jean Dieudonné and Laurent Schwartz took him under their wing. After his initial work in functional analysis, which still remains fundamental, he turned to algebraic geometry with great success, a story that is now well-known.

Grothendieck was one of the first two professors appointed to the Institut des Hautes Études Scientifiques (IHES) in 1959, where he obtained most of his results and published, with the help of Jean Dieudonné, the other professor at the IHES, the famous *EGA* (*Éléments de géométrie algébrique*). He directed the seminar on algebraic geometry, which resulted in a publication of more than 5,000 pages coauthored with some of his students, known as *SGA* (*Séminaire de Géométrie Algébrique du Bois Marie*). Grothendieck was awarded the Fields Medal at the International Congress of Mathematicians in 1966, but did not travel to Moscow to receive it. In 1988, he won the prestigious and well-funded Crafoord Prize, which he refused.

Grothendieck left the IHES in 1970 after discovering that the institute benefited from military funding and launched his own ecological crusade, first through the journal *Survivre*, and then later *Survivre et vivre*. But Grothendieck not only left the IHES, he also left the world of mathematics and, in particular, his students. He returned to mathematics in 1983, but in a very different style, with his publications "*Esquisse d'un programme*" (Sketch of a Programme) and "*À la poursuite des champs*" (Pursuing Stacks).[17] After a year at the Collège de France, he was appointed professor in Montpellier, where he worked until his retirement in 1988. He spent his final years living in the countryside in almost total seclusion, until his death in 2014 at the age of 86.

As we have seen, Grothendieck is the author of a considerable body of mathematical work. But he is also the author of significant literary works. Among them is *R&S*, which was published by Gallimard in January 2022 after having been widely distributed on the internet since Grothendieck first wrote the text in 1986. Amounting to more than 1,900 pages, the book deals with many subjects: the author's journey as a mathematician, his passions, his illusions and disillusions, the process of creation, and a thousand other topics. It also includes long passages on Yin and Yang, feminine and masculine ways of doing mathematics, the mother, the father and child, dreams, and so on. A large part of the text is devoted to a revelation he is said to have experienced in 1976 and a long period of meditation that followed. It is a kind of self-analysis tinged, it has to be said, with a certain degree of paranoia. A recurring theme is the sense of betrayal he felt toward his former students, which is manifested in his work being ignored and forgotten. The words "funeral," "deceased," "hearse," "massacre," and "gravedigger," and so on, quickly become omnipresent after their appearance in the table of contents. More generally, the book denounces a loss of ethics among the entire mathematical community.

Grothendieck explains to the reader that mathematics "was better before"—that is, prior to 1960—as if the older generation was irreproachable! In fact, on the contrary, it can be said that mathematicians have become much more honest since the 1990s. The source of this miracle has a name: *arXiv*. It is now becoming ever more difficult to appropriate the ideas of others, although, of course, it is still possible to some degree. The institution of mathematics itself has also been greatly improved, or at least has been greatly transformed. The system of mandarins that dominated French mathematics until the 1970s, from which Grothendieck did not experience any difficulties and about whom he does not say a word, has practically disappeared.

Grothendieck, who is very self-critical throughout the text, sometimes ponders whether he might have been arrogant or even contemptuous of those around him during his heyday in the 1960s and 1970s. Despite these concerns, it is clear that he cares little about ingratiating himself with his readers. Instead, he offers a book of more than 1,900 pages, while in response to a question about the IHES library in its early days, he remarks: "We don't read books, we write them!"[18] *R&S* contains many contradictions that are only partly corrected by a series of Notes—some of which, despite being of particular importance, are not included in this new edition. Addressing these contradictions properly would undoubtedly have required the text to be completely rewritten.

Grothendieck is not paralyzed by any sense of false modesty:

> The thing that struck me is that I do not remember having known, even from the allusions of friends or colleagues who are better versed in history, of a mathematician apart from myself who contributed such a multiplicity of innovative ideas, not more or less disjointed from one another, but as part of a vast unifying vision (as was the case for Newton and for Einstein in physics and cosmology, and for Darwin and for Pasteur in biology).[19]

Elsewhere, he writes: "It would seem that, as a servant of a vast unifying vision born in me, I am 'one of a kind' in the history of mathematics from its origin to the present day."[20]

Although the writing style is not lacking in inspiration, it is nonetheless uneven and sometimes—deliberately—familiar. Grothendieck is not le duc de Saint-Simon.

*The following analysis will focus only on the content concerning mathematics and the world of mathematicians.*

In the text, Grothendieck complains at length that his ideas have been plundered by his former students without reference to their master or that they have simply been erased and forgotten. These assertions are not always supported by solid arguments or precise references. But, above all, it is the nature of discoveries to be trivialized and their author forgotten, and all the more so when the underlying ideas are often, in hindsight, obvious.

Grothendieck's reproaches are addressed to all his pupils, and particularly to Deligne—whose name is almost always preceded by the words "my friend," insinuating "my former friend"—and to Verdier. It is quite possible to imagine that Deligne was only lightly involved with Grothendieck's authorship of the motives or that the "Verdier duality" already mentioned could just as well be called the "Grothendieck duality." But, otherwise, everyone knows that it was Grothendieck who invented schemes, motives, Grothendieck topologies, topoi, and, above all, that he imposed the functorial point of view via the six operations and the derived categories. Everyone knows that it is thanks to the machinery devised by Grothendieck that Deligne was able to prove André Weil's last conjecture.

In support of his claims about the total loss of ethics in the mathematical community from the 1970s onwards, Grothendieck's entire argument is based on the unique testimony of one and only one mathematician who came to see him several times at his home in the countryside.

It is common practice in ethnology to rely on an informant from the group being studied and who speaks the language. The problem is that the informant may not always be all that reliable and can, in fact, say anything. Here it is an even worse situation, since the informant declares himself to be the first person affected by the story he is going to tell, namely the Riemann–Hilbert (R–H) correspondence.

THE INFORMANT was able to convince Grothendieck that he was, in a certain sense, his spiritual son—"a continuator of my work."[21] Furthermore, the informant persuaded Grothendieck that he had been able to demonstrate the R–H correspondence—without the slightest advice and in the face of indifference, if not outright hostility, from all around him.[22] And that he had done so using the language of derived categories for the first time in this field. Speaking of this informant, Grothendieck also writes that "his pioneering work since 1972 has been done in complete solitude."[23]

*All of this is grossly untrue.*

The informant did his postgraduate thesis in 1974 under my direction and on a subject that I had proposed to him. He had largely benefited from a private talk given by Masaki Kashiwara in 1975 when the informant was starting to prepare his first article, which makes no mention of this decisive talk. He also benefited from a copy of Kashiwara's 1970 thesis[24]—written in Japanese, but there is no lack of translators—and from Christian Houzel's repeated advice throughout his *thèse d'état*. As for the derived categories, they appear on the first page of the foundational article by Mikio Sato, Takahiro Kawai, and Kashiwara published in 1973.[25]

The R–H correspondence is an "equivalence of categories" formulated by Kashiwara in 1975 and was demonstrated by the same author in 1980.[26]

In passing, it is worth noting that interesting equivalences of categories build bridges between different a priori unrelated fields of mathematics: here, the partial differential equations of analysis and the constructible sheaves of algebraic topology. Another more recent and very important equivalence is provided by Maxim Kontsevich's "mirror symmetry," which links complex geometry and symplectic geometry.

Grothendieck's entire statement about the role of his protégé in the story of the R–H correspondence, scattered and repeated throughout the 1,900 pages of his original text, is therefore based on false testimonies. However, with astonishing naïveté, our author takes everything that his interlocutor tells him at face value and he is quoted more than 200 times. Grothendieck goes on a crusade against Kashiwara, whom he goes so far as to call "a ringleader [*caïd* in the original French] from across the Pacific,"[27] even though he had never communicated with Kashiwara and had only a very fragmentary knowledge of his work. By extension, it is the entire Sato school that is labeled "ringleaders from across the Pacific."[28] Without being exhaustive, which would be impossible in any practical sense without copying the entire book, consider the following quote from Note 458, which is not lacking in unintentional irony: "the Sato school is said to have initiated the method of surrounding itself with obscurity in order to dominate."[29]

In 1981, Grothendieck reached the peak of his resentment with an event he referred to as *"le Colloque Pervers"* (the Perverse Colloquium), a historically important conference to which his protégé was not invited. It was on this occasion that Alexander Beilinson, Joseph Bernstein, and Deligne—not counting Ofer Gabber who refused to associate his name with it for obscure reasons—introduced perverse sheaves.[30] The idea of these sheaves—which are not sheaves in the strict sense, but complexes of sheaves, hence, perhaps, the adjective—arises naturally from the R–H correspondence. Their definition already appears implicitly in Kashiwara's 1975 text. Grothendieck was outraged that his protégé was completely ignored at this colloquium when he should have been its star. If no reference was made to the authorship of R–H at the event, it was probably because the mathematical community was aware of the controversies surrounding it and nobody wanted to be involved. But if one reads what Grothendieck writes in the additions to *R&S*, the notable absentee at this colloquium was not, in fact, his informant, but *Kashiwara*! In other words, all the pages and pages of indignation in the original text are either entirely misplaced, or do not defend the right people. There is no doubt in my mind that Sato and his student Kashiwara were unjustly ignored, or maybe even misunderstood, by the French school during the 1980s, and by Bourbaki in particular, but that is another story.[31]

In the edition published by Gallimard, Grothendieck cautiously goes back on his assertions concerning the authorship of R–H and is willing to acknowledge that Kashiwara could have played a role in it,[32] perhaps even a role equivalent to that of his informant. Finally, at the end of Part III, he offers "my most sincere apologies" to Kashiwara.[33] But 1,500 pages later, none of these admissions preclude Grothendieck from insulting Kashiwara once again.

In 1986, I was aware of this part of Grothendieck's text—concerning the ringleaders, or *caïds*, as he referred to them, from the other side of the Pacific—and I wrote to him about it on January 16. An important correspondence followed that continued for several months until around the end of March. Supported by the testimony of Christian Houzel, I believe that I successfully convinced Grothendieck that his version of R–H was completely false.

In a series of additions, adding about twenty pages to the initial text of *R&S*, some extracts from which are presented below,[34] Grothendieck went back completely on what he had written earlier. He finally affirmed that it was indeed Kashiwara who first formulated the R–H correspondence in 1975 and that it was also Kashiwara who gave the first outline of a proof in 1980. It is to Grothendieck's credit that he admits his error of judgement, but because it had been propagated throughout 1,900 pages, it was an error that was not easily rectified. Grothendieck chose to address the problem using additions and footnotes, but unfortunately, apart from a flat apology,[35] none of these amendments appear in the published book.

Consider, for example, the following additions.

Grothendieck, writing on May 9, 1986:

After the provisional distribution of *Récoltes et semailles*, from October last year, I was contacted by Pierre Schapira, and then by Christian Houzel, to point out some glaring inaccuracies in the version of events presented in *Récoltes et semailles*. The situation was clarified considerably during correspondence with both of them, which continued between January and March of this year. It now appears to me that in the "[Zoghman] Mebkhout version" (which was not lacking in internal consistency) the true, the tendentious and the downright false are inextricably mixed.

Grothendieck, also dated May 15, 1986:

In retrospect, I am convinced that Kashiwara cannot be reproached for the slightest incorrectness in this case. In his presentation, he gives a statement and a first sketch of a proof of a theorem, which he had indeed been the first to conjecture as early as 1975... Moreover, he has the correction to specify, as early as page 2: "Let us note that the theorem is also proved by Mebkhout by a different way." This was even "lending to the rich," because the previous month, in his note to the CRAS [*Comptes Rendus de l'Académie des Sciences*] of 3 March 1980, Mebkhout had expressed himself in the hypothetical form "we hope to show that...," and without making the slightest allusion to it...

All of this raises some questions about the editorial process that led to the publication of *R&S* by Gallimard in January 2022. In the foreword, which is dated January 1986,[36] Grothendieck thanks Christian Bourgois and Stéphane Deligeorges for including his text in the *Épistémè* collection. What happened between this date and the publication by Gallimard? And, above all, why do Grothendieck's additions—incorporated prior to May 29, 1986—not appear in the final published version, while the brief apology that does appear clearly proves that the Gallimard edition includes other elements added after January 1986?[37]

In this review, I have focused on the mathematical sections of the book and the passages that discuss the history of the R–H correspondence. The latter is far from being an anecdotal component in the book. Indeed, Grothendieck refers to it constantly—it is a leitmotif. Unfortunately, and as he admits with great frankness, Grothendieck was misled by an informant lacking in objectivity, to say the absolute least. With a disarming and, in a certain sense, admirable degree of naivety, he also admits that he never imagined that the information he was being fed could be biased or incomplete, let alone downright false. Between 1955 and 1970, Grothendieck lived in a world of pure ideas. He was immersed in mathematics to an extent that is hard to imagine. When he emerged from the noosphere into the real world—that is, the social world—one can only imagine the harsh shock of everyday life and how crushed he felt by what he perceived as a loss of ethics in the world of science. But why should this world be any different from the rest of society? The rigor of science has never been reflected in its practitioners—examples are legion.

Science is a great devourer of men and characters.[38]

1. The author would like to thank Leila Schneps for her critical and constructive advice.
2. Alexander Grothendieck, *Récoltes et Semailles: I, II. Réflexions et témoignage sur un passé de mathématicien* (Paris: Gallimard, 2022), 439.
3. See the historical article by Christian Houzel "Les débuts de la théorie des faisceaux," in Masaki Kashiwara and Pierre Schapira, *Sheaves on Manifolds, Grundlehren der Mathematischen Wissenschaften, vol. 292* (Berlin: Springer-Verlag, 1990), doi:10.1007/978-3-662-02661-8.
4. Roger Godement, *Théorie des faisceaux* (Paris: Hermann, 1958).
5. Jean-Pierre Serre, "Faisceaux Algebriques Coherents," *Annals of Mathematics, 2nd Series* 61, no. 2. (1955): 197–278, doi:10.2307/1969915 .
6. Alexander Grothendieck, "Sur quelques points d'algèbre homologique," *Tōhoku Mathematical Journal* 9, no. 3 (1957): 119–21, doi:10.2748/tmj/1178244774. This article is analyzed in detail in Rick Jardine, "Tōhoku," *Inference* 1, no. 3 (2015), doi:10.37282/991819.15.13.
7. Mike Artin, Alexandre Grothendieck, and Jean-Louis Verdier, *Théorie des topos et cohomologie étale des schémas, Lecture Notes in Mathematics*, vols. 269, 270, 305 (Berlin: Springer-Verlag, 1972–73).
8. Ralf Krömer, "La « machine de Grothendieck » se fonde-t-elle seulement sur des vocables métamathématiques ? Bourbaki et les catégories au cours des années cinquante," *Revue d'histoire des mathématiques* 12 (2006): 119–62.
9. Samuel Eilenberg and Saunders Mac Lane, "Natural Isomorphisms in Group Theory," *Proceedings of the National Academy of Sciences* 28 (1942): 537–43, doi:10.1073/pnas.28.12.537; Samuel Eilenberg and Saunders Mac Lane, "General Theory of Natural Equivalences," *Transactions of the American Mathematical Society* 58 (1945): 231–94, doi:10.1090/S0002-9947-1945-0013131-6.
10. Henri Cartan and Samuel Eilenberg, *Homological Algebra* (Princeton, NJ : Princeton University Press, 1956).
11. Grothendieck, *Récoltes et Semailles*, 158.
12. Grothendieck, *Récoltes et Semailles*, 42, 377.
13. Alexander Grothendieck, "A la poursuite des champs," (1987).
14. Pierre Cartier, "A Country Known Only by Name," *Inference* 1, no. 1 (2014), doi:10.37282/991819.14.7. For additional details, see "Sascha Shapiro," *Wikipedia*.
15. Allyn Jackson, "*Comme Appelé du Néant*—As if Summoned from the Void: The Life of Alexandre Grothendieck," *Notices of the AMS* 51, no. 4 and 51, no. 10 (2004): 1,038–54 and 1,196–212; Winfried Scharlau, "Who is Alexander Grothendiek?," *Notices of the AMS* 55, no. 8 (2008): 930–41; Leila Schneps, "The Grothendieck Circle" (grothendieckcircle.org).
16. The wording "extreme-left militant in Germany in the 1920s" does not really have the same meaning as its transposition a century later.
17. Grothendieck, "A la poursuite des champs."
18. Jackson, "*Comme Appelé du Néant*," 1,050.
19. Grothendieck, *Récoltes et Semailles*, 935.
20. Grothendieck, *Récoltes et Semailles*, 94.
21. Grothendieck, *Récoltes et Semailles*, 1,664.
22. Grothendieck, *Récoltes et Semailles*, 413.
23. Grothendieck, *Récoltes et Semailles*, 1,663.
24. Masaki Kashiwara, "Algebraic Study of Systems of Partial Differential Equations," Master's thesis (Tokyo University, 1970), *Mémoires de la Société Mathématique de France* 63 (1995).
25. Mikio Sato, Takahiro Kawai, and Masaki Kashiwara, "Microfunctions and Pseudo-Differential Equations, Hyperfunctions and Pseudo-Differential Equations," in *Proceedings of a Conference at Katata, 1971; Dedicated to the Memory of André Martineau, Lecture Notes in Mathematics*, vol. 287 (Berlin: Springer-Verlag, 1973), 265–529.
26. A few words on the R–H correspondence follow. The modern formulation uses the theory of *D*-modules, a theory that was intuited by Sato in the 1960s, and fully implemented by Kashiwara in his thesis. (A related theory was developed independently by Joseph Bernstein.) In everyday language, a coherent *D*-module means a system of partial differential equations with holomorphic coefficients. Holonomic modules are the > 1 dimensional version of ordinary differential equations, and among them regular holonomic modules generalize the classical notion of Fuchsian equations. In 1975, Kashiwara showed that the functor Sol, which associates the complex of its holomorphic solutions to a holonomic module, takes its values in constructible sheaves, the sheaves which behave locally as a direct sum of constant sheaves along a stratification. In the same year he conjectured that there exists a triangulated subcategory of the holonomic modules, the regular holonomic modules, on which the functor Sol induces an "equivalence of categories." Kashiwara proved his conjecture in 1980 and gave a detailed account of the main steps in his proof at the École Polytechnique seminar, which was published. His proof uses Heisuke Hironaka's singularity resolution theorem and the precursor work of Deligne who treated the special case of "meromorphic connections."
27. Grothendieck, *Récoltes et Semailles*, 1,656.
28. Grothendieck, *Récoltes et Semailles*, 1,651.
29. Grothendieck, *Récoltes et Semailles*, 1,650.

30. Alexander Beilinson, Joseph Bernstein, and Pierre Deligne, "Faisceaux pervers, Analysis and topology on singular spaces, I" (Luminy, 1981), *Astérisque*, no. 100 (Paris: Societé Mathematique de France, Paris, 1982): 5–171.
31. Pierre Schapira, "Mikio Sato, a Visionary of Mathematics," *Notices of the AMS* 54, no. 2 (2007); Pierre Schapira, "Fifty Years of Mathematics with Masaki Kashiwara," Proceedings of the International Congress of Mathematicians, Rio de Janeiro, 2018, vol. 1, Plenary Lectures (Singapore: World Scientific, 2018).
32. Grothendieck, *Récoltes et Semailles*, 163.
33. Grothendieck, *Récoltes et Semailles*, 164.
34. These additions have just been posted on Leila Schneps's website, "The Grothendieck Circle" (grothendieckcircle.org).
35. Grothendieck, *Récoltes et Semailles*, 163.
36. Grothendieck, *Récoltes et Semailles*, 9–15.
37. On the recently updated "Grothendieck Circle" (grothendieckcircle.org) website, it is written that Gallimard editions will soon publish a new version of *R&S* augmented with these Grothendieck additions.
38. Adapted freely from a famous quotation by Leon Trotsky.

Sorbonne Université, Cnrs IMJ-PRG
pierre.schapira@imj-prg.fr
http://webusers.imj-prg.fr/~pierre.schapira

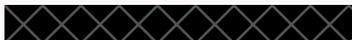
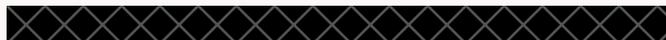